\documentclass[letterpaper, 10 pt, conference]{ieeeconf}  % Comment this line out if you need a4paper
\usepackage{cite, graphicx, amsmath, amssymb, mathrsfs, float, subfig,color}
\newtheorem{theorem}{Theorem}
\renewcommand{\baselinestretch}{0.96}

\IEEEoverridecommandlockouts                              % This command is only needed if 
\title{\LARGE \bf
Robust Decentralized Voltage Control of DC-DC Converters with Applications to Power Sharing and Ripple Sharing}

\author{Mayank Baranwal$^{1,a}$, Srinivasa M. Salapaka$^{1,b}$ and Murti V. Salapaka$^{2,c}$% <-this % stops a space
\thanks{$^{1}$Department of Mechanical Science and Engineering, University of Illinois at Urbana-Champaign, IL 61801, USA}%
\thanks{$^{2}$Department of Electrical and Computer Engineering, University of Minnesota, Minneapolis, MN 55455 USA}%
\thanks{$^{a}${\tt\small baranwa2@illinois.edu}, $^{b}${\tt\small salapaka@illinois.edu}}
\thanks{$^{c}${\tt\small murtis@umn.edu}}        
}

\begin{document}

\maketitle
\thispagestyle{empty}
\pagestyle{empty}

%%%%%%%%%%%%%%%%%%%%%%%%%%%%%%%%%%%%%%%%%%%%%%%%%%%%%%%%%%%%%%%%%%%%%%%%%%%%%%%%
\begin{abstract}
This paper addresses the problem of output voltage regulation for multiple  DC-DC converters connected to a grid, and prescribes a robust scheme for sharing power among different sources. Also it develops a method for sharing 120 Hz ripple among  DC power sources in a prescribed proportion, which accommodates the different capabilities of DC power sources to sustain the ripple. We present a decentralized control architecture, where a  nested  (inner-outer) control design is used at every converter. An interesting aspect of the proposed design is that the analysis and design of the entire multi-converter system can be done using an equivalent single converter system, where  the multi-converter system inherits the  performance and robustness achieved by a design for the single-converter system. Another key aspect of this work is that the voltage regulation problem is addressed as a disturbance-rejection problem, where {\em unknown} load current is viewed as an external signal, and thus, no prior information is required on the nominal loading conditions. The control design is obtained using robust optimal-control framework. Case studies presented show the enhanced performance of prescribed optimal controllers.
\end{abstract}

%%%%%%%%%%%%%%%%%%%%%%%%%%%%%%%%%%%%%%%%%%%%%%%%%%%%%%%%%%%%%%%%%%%%%%%%%%%%%%%%
\section{INTRODUCTION}

\begin{figure}[!t]
	\centering
	\includegraphics[width=\columnwidth]{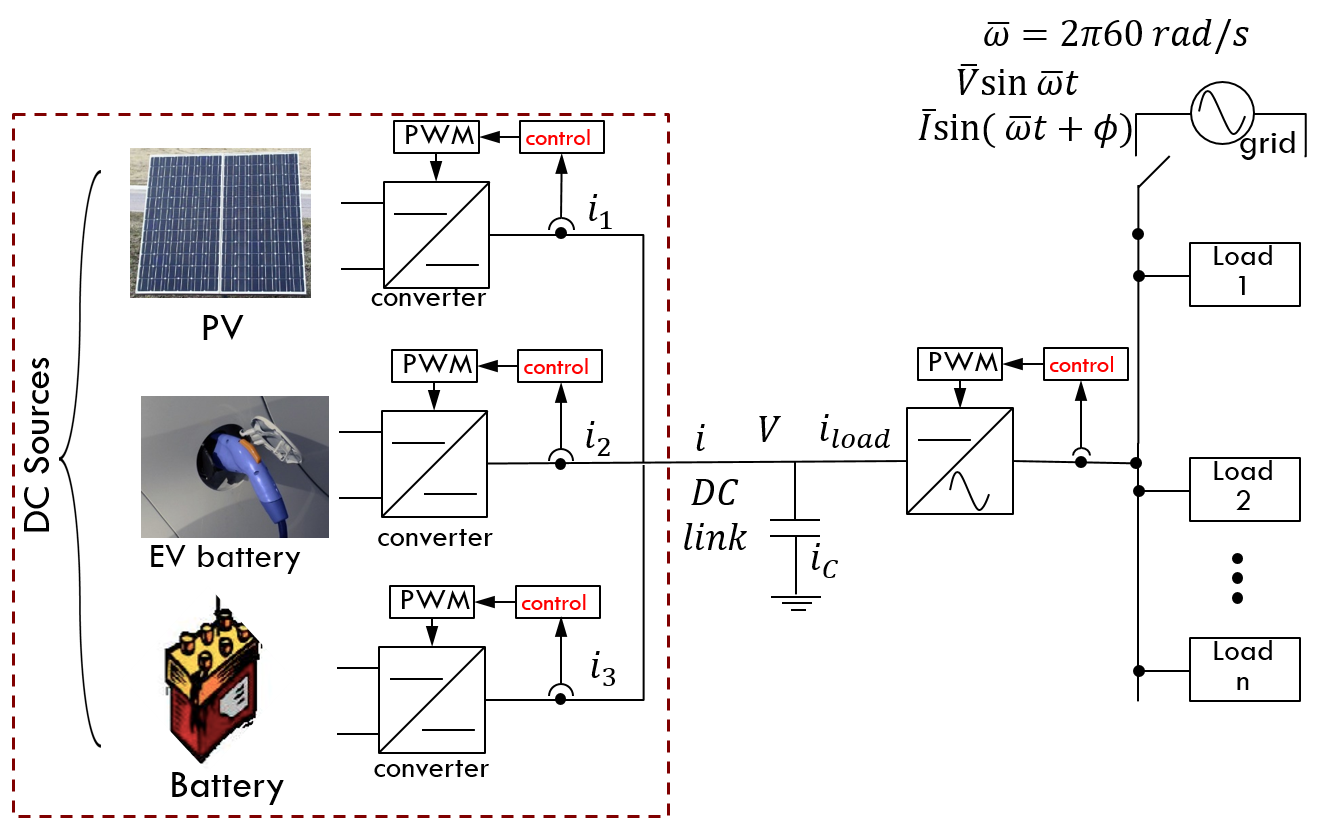}
	\vspace{-0.5em}
	\caption{{\small A schematic of a microgrid. An array of DC sources provide power for AC loads. Power sources  provide  power at DC-link, their common output bus, at a voltage that is regulated to a set-point. The control system at the respective DC-DC converter that interfaces with a source is responsible for regulating the voltage at he DC-link.    An inverter that connects to the DC-link converts the total current from the sources at the regulated voltage to  alternating current (AC)   at its output to satisfy the power demands of the AC loads. This paper describes an approach for control design of the multiple converters systems associated with power transfer from sources to the DC-link (shown by the dotted line).}} 
	\label{fig:MCS}
	\vspace{-1.5em}
\end{figure}

In power network topologies, especially in microgrids \cite{hatziargyriou2007microgrids}, multiple DC power sources connected in parallel (see Figure \ref{fig:MCS}), each interfaced with DC-DC converter, provide power at their common output, the DC-link, at a regulated voltage; this power can directly feed DC loads or be used by an inverter to interface with AC loads . Voltage controllers form an integral component of DC-DC converters in such systems. A paralleled architecture for multiple power sources  is preferred  since it enables higher output power, higher reliability and ease of use \cite{thottuvelil1997analysis}. Here two main control architectures are adopted - (1) {\em master-slave} control, where the voltage regulation error from the master converter is utilized to provide an error signal to all the parallel connected converters \cite{thottuvelil1997analysis, sha2010cross}, (2) {\em decentralized} control, where each converter utilizes an independent and variable voltage reference depending on the output of each unit \cite{panov1997analysis, guerrero2007decentralized}. Irrespective of the control framework, controllers at each converter are to be designed such that the  voltage  at the DC-link is regulated at a prescribed setpoint. Another important control objective is to ensure that the DC sources provide power in a prescribed proportion, which may be dictated by their power ratings or external economic criteria.  The main challenges arise from the uncertainties in the size and the schedules of loads, the complexity of a coupled multi-converter network, the uncertainties in the model parameters at each converter, and the adverse effects of interfacing DC power sources with AC loads, such as the $120$ Hz ripple that has to be provided by the DC sources. 
 
Problems related to robust and optimal design of converter controllers have received recent attention. In \cite{olalla2010lmi}, a linear-matrix-inequality (LMI) based robust control design for boost converters has resulted in significant improvements over conventional PID-based controllers. Use of $\mathcal{H}_\infty$ framework in context of inverter systems has also been studied in \cite{weiss2004h, hornik2011current, salapaka2014viability}. While the issue of current sharing is extensively studied \cite{panov1997analysis, ashley1996current, batarseh2010method, wu1993load}, most methods assume a single power source. A systematic control design that addresses all the challenges and objectives for the multi-converter control is still lacking.

In this paper, we develop a control architecture that addresses the following primary objectives for multi-converter control - (1) {voltage regulation} at the DC-link while guaranteeing { robustness} of the closed-loop system  performance to load and parametric uncertainties, (2) { prescribed power sharing} among a number of parallel converters, (3) { controlling the tradeoff between $120$ Hz ripple on  the total current provided by the power sources and the ripple on the DC-link voltage}, and (4) $120$ Hz DC {output ripple sharing}  among converters.  The tradeoffs with $120$ Hz ripple in objectives (3) and (4) result from a direct consequence of interfacing DC-power sources with AC loads (see Figure \ref{fig:MCS}). If we assume negligible power losses at the inverter and the load bus,  the power provided by power sources at the DC-link should equal to the power consumed by the AC loads, that is, $Vi=(\bar V\sin\bar\omega t)(\bar I\sin(\bar\omega t+\phi))=\dfrac{\bar{V}\bar{I}}{2}(\cos\phi+\cos(2\bar\omega t+\phi)$. Since the instantaneous power has a $2\bar\omega=120$ Hz ripple, the current $i=\sum i_k$ at the DC-link has to provide for this $120$ Hz component. The total ripple demand posed by the AC-grid side  is met partly by the ripple current sourced by a capacitor $i_C$  which reduces the magnitude of the ripple current to be provided by the DC source via the converters. However, greater the ripple magnitude in the capacitor current greater will be the ripple in the capacitor voltage, thus adversely affecting  voltage regulation. Therefore a compromise  has to be reached in the allowable ripple in the capacitor current and the ripple provided by the DC sources. The control scheme presented in the article provides a ``knob" to adjust the relative ratio of how the $120$ Hz ripple is shared between the two quantities - the sourced current $i$ and the output voltage $V$. Moreover in a scenario where multiple and different types of DC sources are employed, it is often the case that the tolerance to ripple varies. Here it becomes important to allocate greater percentage of ripple load to tolerant DC sources while reducing the ripple load on vulnerable DC sources. The article presents a controller synthesis procedure where the $120$ Hz ripple on the current $i$ can be shared among the paralleled DC sources ($i_k$) in a pre-specified proportion.

An important aspect of the proposed control architecture is that it is decentralized and addresses all the objectives {\em simultaneously}. Moreover we show that for the control approach described in the paper, the control design and the closed-loop analysis of the multi-converter system can be completely characterized in terms of an appropriate single-converter system; thereby significantly reducing the complexity in addressing multi-converter systems. This architecture exploits structural features of the paralleled multi-converter system, which results in a modular and yet coordinated control design. For instance, it exploits that  the voltage regulation objective is common to all the converters, and that the differences in demands on different converters are mainly in terms of their output currents; accordingly   at each converter, it 
 employs a nested  (outer-voltage inner-current) control structure  \cite{erickson2007fundamentals}, where the outer loop is responsible for robust voltage regulation and the inner loop for shaping the currents. The structure of control  for each inner-loop is so chosen that the  entire closed-loop multi converter system  can be reduced to an equivalent  single-converter system  in terms of the transfer function  from the desired regulation setpoint $V_{des}$ to the voltage $V$.  Furthermore, for the outer-control design, the load current is treated as an external disturbance and the voltage regulation problem is cast as a disturbance rejection problem in an optimal control setting. This design, besides achieving the voltage-regulation objective, provides robustness to deviations from the structural assumptions in the control design. Note that this viewpoint is in contrast to typical methods in  existing literature, where voltage regulation in presence of { unknown} loads is addressed either using { adaptive control}  \cite{oucheriah2013pwm}, or by letting the voltage droop in a controlled manner.

\section{Modeling of converters}\label{sec:avg_model}
In this section, we provide dynamic models for DC-DC converters, which convert a source of direct current (DC) from one voltage level to another. The  models presented below depict dynamics for signals that are averaged over a switch cycle.

\begin{figure*}[!t]
	\begin{center}
	\begin{tabular}{ccc}
	\includegraphics[width=0.64\columnwidth]{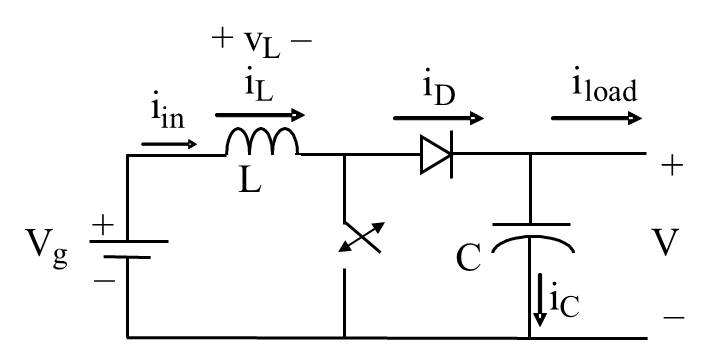}&\includegraphics[width=0.64\columnwidth]{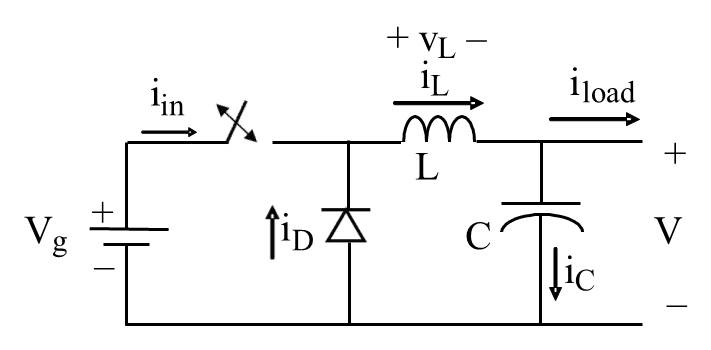}&\includegraphics[width=0.64\columnwidth]{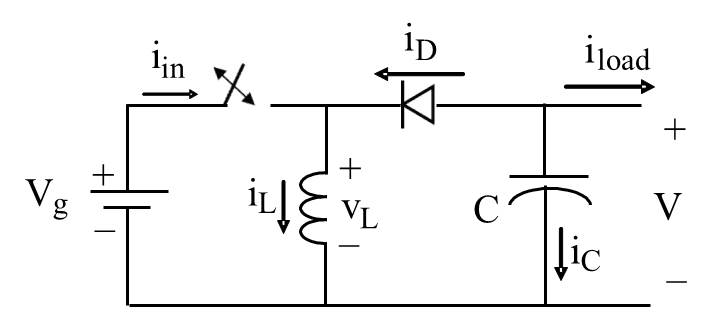}\cr
	(a)&(b)&(c)\end{tabular}
	\vspace{-0.5em}
	\caption{{\small Circuit representing (a) Boost converter, (b)Buck converter, and (c) Buck-Boost converter. Note that $i_{\text{load}}$ includes both the nominal load current, as well as ripple current. The converters are assumed to operate in continuous-conduction-mode (CCM). Boost converters step up the voltage at the output, while buck converters step down the voltage. A buck-boost converter can achieve both the objectives.}}
	\label{fig:model_boost}
	\vspace{-2em}
	\end{center}
\end{figure*}
A schematic of the Boost converter is shown in Fig. \ref{fig:model_boost}(a). %Under assumed to be operating in CCM and semiconductors switches are assumed to be ideal.   
 A dynamic model (averaged over switching cycles) is given by,
\begin{small}
\begin{equation*}\label{eq:boost_model1}
	L \dot{i_L}(t) = -(1-d(t))V(t) + V_g, \quad	C \dot{V}(t) = (1-d(t))i_L(t) - i_{\text{load}},
\end{equation*}\end{small} where $d(t)$ represents the duty-cycle (or the proportion of {\em ON} duration) at time {\em t}, which by defining $d'(t) = 1-d(t)$ and $D' = \frac{V_g}{V_{\text{des}}}$ can be rewritten as
\begin{small}
\begin{equation*}\label{eq:boost_model2}
	L\dot{i_L}(t) = \underset{\tilde{u}(t) := V_g - u(t)}{\underbrace{-d'(t)V(t) + V_g}}, \quad
	C\dot{V}(t) = \underset{\approx D'}{\underbrace{(D'+\hat{d'})}}i_L(t) - i_{\text{load}}.
\end{equation*}\end{small} Here $V_{\text{des}}$ represents the desired output voltage and $\hat d=d'(t)-D'$ is typically very small, which allows for a 
 linear approximation around the nominal duty-cycle, $D = 1-D'$ given  by,
\begin{small}
\begin{equation*}\label{eq:boost_model3}
	L\frac{di_L(t)}{dt} = \tilde{u}(t), \quad
	C\frac{dV(t)}{dt} \approx D'i_L(t) - i_{\text{load}}.
\end{equation*}
\end{small}

Fig. \ref{fig:model_boost}(b) depicts the circuit schematic of a buck converter with an ideal switch. The averaged model of a buck converter is given by,
\begin{small}
\begin{equation*}\label{eq:buck_model}
	L\frac{di_L(t)}{dt} = \underset{\tilde{u}(t) := -V(t) + u(t)}{\underbrace{-V(t) + d(t)V_g}}, \quad
	C\frac{dV(t)}{dt} = i_L(t) - i_{\text{load}}.
\end{equation*}
\end{small}

The electronic circuit of a buck-boost converter is shown in Fig. \ref{fig:model_boost}(c). As in case of a boost converter, we define nominal duty-cycle, $D = \frac{V_{\text{des}}}{V_{\text{des}}-V_g} = 1-D'$. A linear approximation of the above dynamical equations yields,
\begin{small}
\begin{equation*}\label{eq:buck_boost_model2}
	L\dot{i_L}(t) = \underset{\tilde{u}(t) := V(t) + u(t)}{\underbrace{V(t) + d(t)(V_g-V(t))}}, \quad
	C\dot{V}(t) \approx -D'i_L(t) - i_{\text{load}}
\end{equation*}
\vspace{-1em}
\end{small}

\section{Control framework for a single-converter}\label{sec:one_control}
\begin{figure}[tphb]
	\centering
	\vspace{-0.5em}
	\includegraphics[width=3.0in]{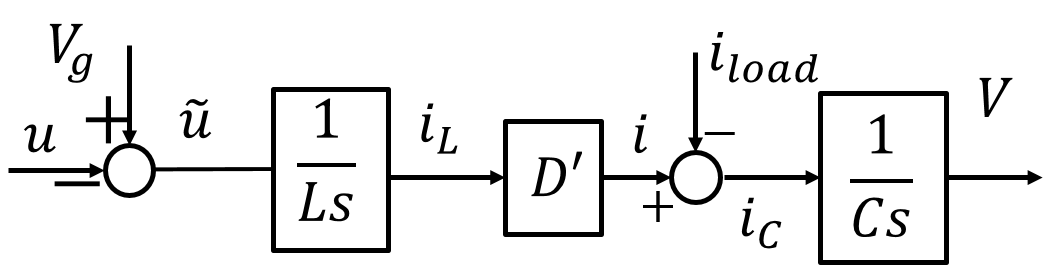}
	\caption{{\small Block diagram representation of a boost-type converter. The control signal $\tilde{u}$ is converted to an equivalent PWM signal to command the gate of the transistor acting as a switch.}}
	\label{fig:plant_representation}
	\vspace{-1em}
\end{figure}

In this section, we describe the inner-current outer-voltage control architecture for a single DC-DC converter. The key objectives of the  control design are - (1) voltage regulation in presence of uncertain loads, and (2) $120$ Hz ripple sharing control between $i_L$ and $i_C$. We first consider the case of a boost converter control design, the dynamics of which is given by 
\begin{small}
\begin{equation}
	i_L(s) = \frac{1}{sL}(V_g(s) - u(s)), \quad
	V(s) = \frac{1}{sC}(D'i_L(s) - i_{\text{load}}(s)),
\end{equation}
\end{small}
 and the corresponding  block-diagram representation of above set of equations is shown in Fig. \ref{fig:plant_representation}, the control objectives are to design $u$ (equivalently $\tilde u$) such that voltage regulation error $V_{\text{des}}-V$ is made small irrespective of load disturbances $i_{\text{load}}$ and  variations in parameters $L$ and $C$, and achieve a prescribed tradeoff between $|i_L(j2\pi120)|$ and $|i_C(j2\pi120)|$.  

These two objectives are achieved using a nested  inner-current outer-voltage control architecture, shown in Fig. \ref{fig:single_inner_controller} (here $G_c = \frac{1}{sL}$ and $G_v = \frac{1}{sC}$). The voltage controller $K_v$ generates a current reference for the current controller $K_c$. The current controller $K_c$ is designed to achieve a high closed inner-loop bandwidth with ripple control as an objective, whereas the voltage controller $K_v$ is designed to achieve a relatively lower closed outer-loop bandwidth  with DC (zero frequency) voltage regulation as its primary objective. We assume that the quantities - output voltage $V$ and inductor current $i_L$ are available for measurement. 

\subsubsection*{Design of the outer-loop controller}\label{sec:robust}
\begin{figure}[tphb]
	\centering
	\includegraphics[width=\columnwidth]{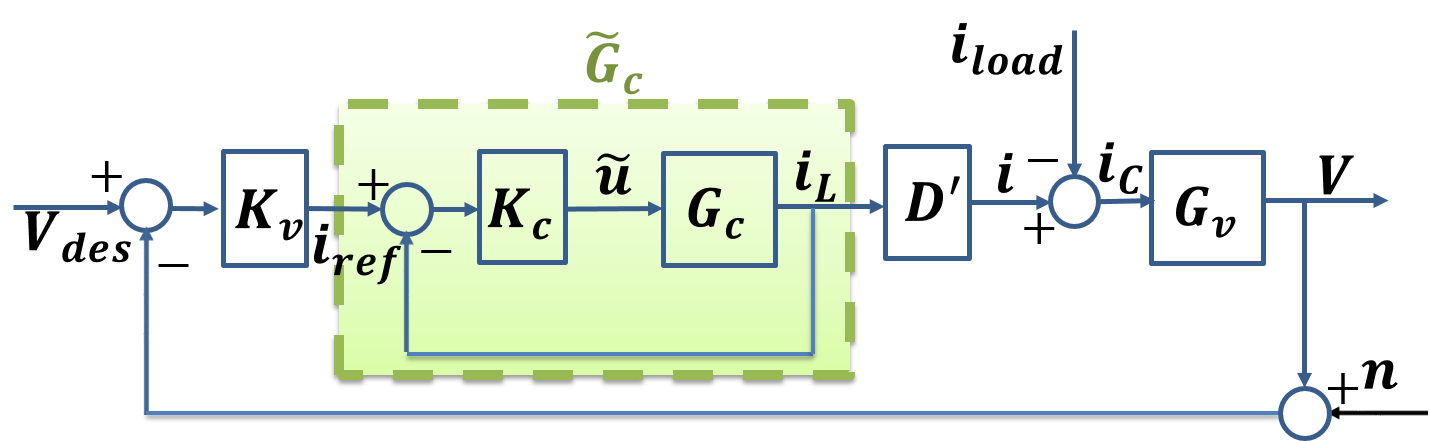}
	\vspace{-2em}
	\caption{{\small Block diagram representation of the inner-outer control design. Exogenous signals $V_{des}$ and $i_{\text{load}}$ represent the desired output voltage and disturbance, respectively. The quantities $V$ and $i_L$ represent the available measurements.}}
	\vspace{-.5em}
	\label{fig:single_inner_controller}
\end{figure}

 For a given controller $K_c$ for the inner-loop, the closed outer-loop signals of interest are given by (see Figure \ref{fig:single_inner_controller}) 
\begin{small}
\begin{eqnarray}
	V_{\text{des}}-V &=& SV_{des} + G_vSd + Tn \label{eq:trade_off1}\\
	i_L &=& \tilde{G}_cK_vSV_{des} + \frac{1}{D'}Td-K_v\tilde{G}_c Sn \label{eq:trade_off2}\\
	i_{\text{ref}} &=& K_v(SV_{des} + G_vSd) - K_vSn, \label{eq:trade_off3}
\end{eqnarray}
\end{small}
where $d$ denotes the load current $i_{\text{load}}$ (shown as disturbance to the plant), $n$ denotes the voltage measurement noise, $\tilde{G}_c$  represents the closed inner-loop transfer function from $i_{\text{ref}}$ to $i_L$,  
and $T(s)$ and $S(s)$ are closed-loop {\em complementary sensitivity} and {\em sensitivity} transfer functions respectively, described by,
\begin{small}
\begin{eqnarray}\label{eq:SnT}
	T(s) &=& (I+G_vD'\tilde{G}_cK_v)^{-1}(G_vD'\tilde{G}_cK_v), \nonumber \\
	S(s) &=& (I+G_vD'\tilde{G}_cK_v)^{-1}. %\nonumber %\\
	%S(s) + T(s) &=& 1,
\end{eqnarray}
\end{small}
The voltage-regulation objective, as evident from  Eq. \ref{eq:trade_off1}, requires designing $K_v$ such that  $|S(j\omega)|$ is small at the frequencies where disturbance $d$ is prominent. However this implies that   the effect of $d$ on the inductor current $i_L$ (see Eq. \ref{eq:trade_off2}) is larger since $|T(j\omega)|$ is larger in those frequencies (from Eq. \ref{eq:SnT}). Therefore there is a fundamental trade-off between voltage regulation and minimizing effects of disturbances (or load current $i_{\text{load}}$) in $i_L$. Also to diminish the effect of noise on voltage regulation, the control design should be such that the closed-loop map $T$ rolls of at frequencies beyond the disturbance bandwidth. Furthermore low $i_{\text{ref}}$ is ensured if $K_vS$ can be made small. The controller $K_v$ is obtained by casting these multiple objectives in the following optimal control problem     \cite{skogestad2007multivariable},
\begin{small}
\begin{equation}\label{eq:stacked_h_infinity}
	\min_{K_v \in \mathscr{K}} \left\| \begin{array}{c}
	W_sS\\
	W_uK_vS\\
	W_tT
	\end{array} \right\|_\infty,
\end{equation}
\end{small}
where the weights $W_{s}$, $W_{t}$ and $W_{u}$ are chosen to reflect the design specifications of robustness to disturbances and parametric uncertainties, tracking bandwidth, and saturation limits on the control signal. For example, the weight function $W_s(j\omega)$ is chosen to be large in frequency range $[0,\omega_{BW}]$ to ensure a small tracking error $e = V_{des}-V$ in this frequency range. The weight function $W_t(j\omega)$ is designed as a high-pass filter to ensure that $T(j\omega)$ is small at high frequencies  to provide  mitigation to measurement noise. The design of constant $W_u$ entails ensuring that the control effort lies within saturation limits. The resulting controller is robust to disturbances up to $\omega_{BW}$, which accounts for variations in load disturbances as well as parametric uncertainties. 

\subsubsection*{Design of the inner-loop controller}\label{sec:inner_controller}
\begin{figure}[tphb]
	\centering
	\includegraphics[width=0.9\columnwidth]{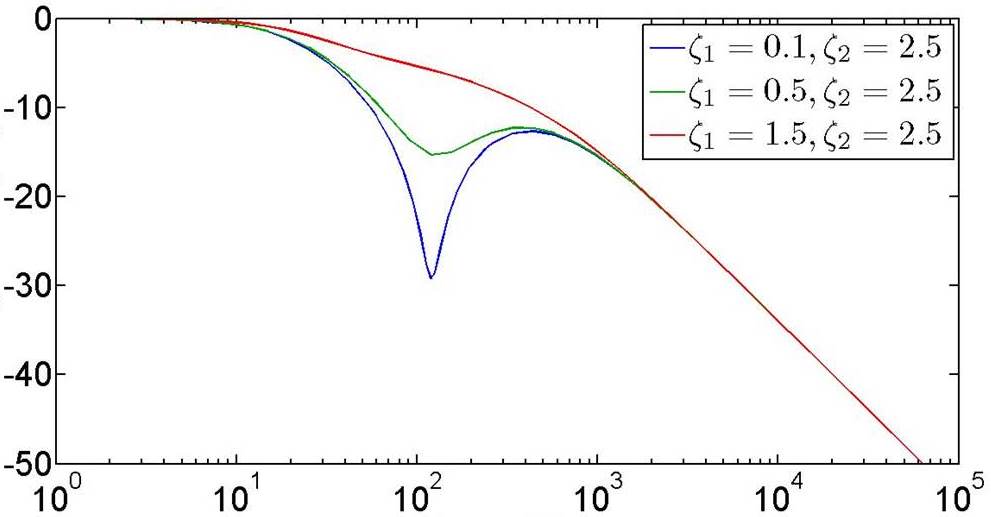}
	\vspace{-0.5em}
	\caption{{\small Bode magnitude plots of the closed-loop plant $\tilde{G}_c$ for various $\zeta_1$ values. $\tilde{\omega}$ is chosen to be 600$\pi$ rad/s. Note that a relatively larger value of $\tilde{\omega}$ is in accordance with choosing a fast inner-current controller.}}
	\label{fig:bode_zeta}
	\vspace{-.5em}
\end{figure}

The outer-loop control design assumed the inner closed-loop $\tilde G_c$. Here we propose an inner-loop control design that results in a second-order transfer function $\tilde G_c$, thereby ensuring a relatively low-order optimal controller $K_v$. The main   objective for designing the inner-loop controller $K_c$ is to decide the trade-off between the $120$ Hz ripple on the voltage $V$ (equivalently on the capacitor current $i_c$) and the output current $i$ (equivalently $i_L$) of the converter.  Accordingly, we design   $K_c$  such that 
\begin{small}
\begin{equation}\label{eq:Gc_tilde}
	\tilde{G}_c(s) = \left(\frac{\tilde{\omega}}{s+\tilde{\omega}}\right)\left(\frac{s^2+2\zeta_1\omega_0 s+\omega_0^2}{s^2+2\zeta_2\omega_0 s+\omega_0^2}\right),
\end{equation}
\end{small}
where $\omega_0 = 2\pi 120$ rad/s. $\tilde{\omega}, \zeta_1, \zeta_2$ are design parameters. Here the parameter $\tilde\omega>\omega_0$ is simply chosen to implement a low-pass filter that attenuates undesirable frequency content in $i_L$ beyond $\tilde\omega$. Note that in this design of $\tilde{G}_c$, there is a notch at $\omega_0=120$ Hz,  the size of this notch is determined by the ratio $\frac{\zeta_1}{\zeta_2}$ (see Figure \ref{fig:bode_zeta}). Lower values of this ratio correspond to larger notches, which  in turn imply smaller $120$ Hz component in $i_L$, since  $\tilde G_c$ represents the inner closed-loop transfer function from $i_{\text{ref}}$ to $i_L$. Furthermore since $i_C=i_{\text{load}}-i_L$,  this in turn implies higher ripples in $i_C$. Thus the ratio $\frac{\zeta_1}{\zeta_2}$ can be appropriately designed to achieve a specified  tradeoff between $120$ Hz ripple on $i_C$ and $i_L$. 
The stabilizing second-order controller $K_c$ that yields the above closed-loop plant $\tilde{G}_c$ is explicitly given by,
\begin{small}
\begin{equation}\label{sec:Kc}
	K_c = L\tilde{\omega}\frac{\left(s^2+2\zeta_1\omega_0 s + \omega_0^2\right)}{\left(s^2+2\zeta_2\omega_0 s + 2(\zeta_2-\zeta_1)\omega_0\tilde{\omega} +\omega_0^2\right)},
\end{equation}
\end{small}
which is again a low-order (second-order) controller design.  

\subsubsection*{Extension to buck and buck-boost converters}\label{sec:ext_to_buck}
The extension of the proposed control design  to Buck and Buck-Boost DC-DC converters is  easily explained after noting that their averaged models are structurally identical to Boost converters, except that the dependence of duty cycles on the control signal $u$ or constant parameter $D'$ are different. The differences in how duty cycles depend on $u(t)$ do not matter from the control design viewpoint since duty cycles for pulse-width modulation are obtained only after obtaining the control designs (that use the averaged models).

\section{Extension to a system of parallel converters}\label{sec:many_control}
In this section we develop a decentralized control framework that achieves voltage regulation, power-sharing, and ripple-sharing among a system of parallel boost converters sharing a common load.

\subsection{Control framework for a system of parallel boost converters}\label{sec:control_scheme}
\begin{figure}[tphb]
	\centering
	\includegraphics[width=3.0in]{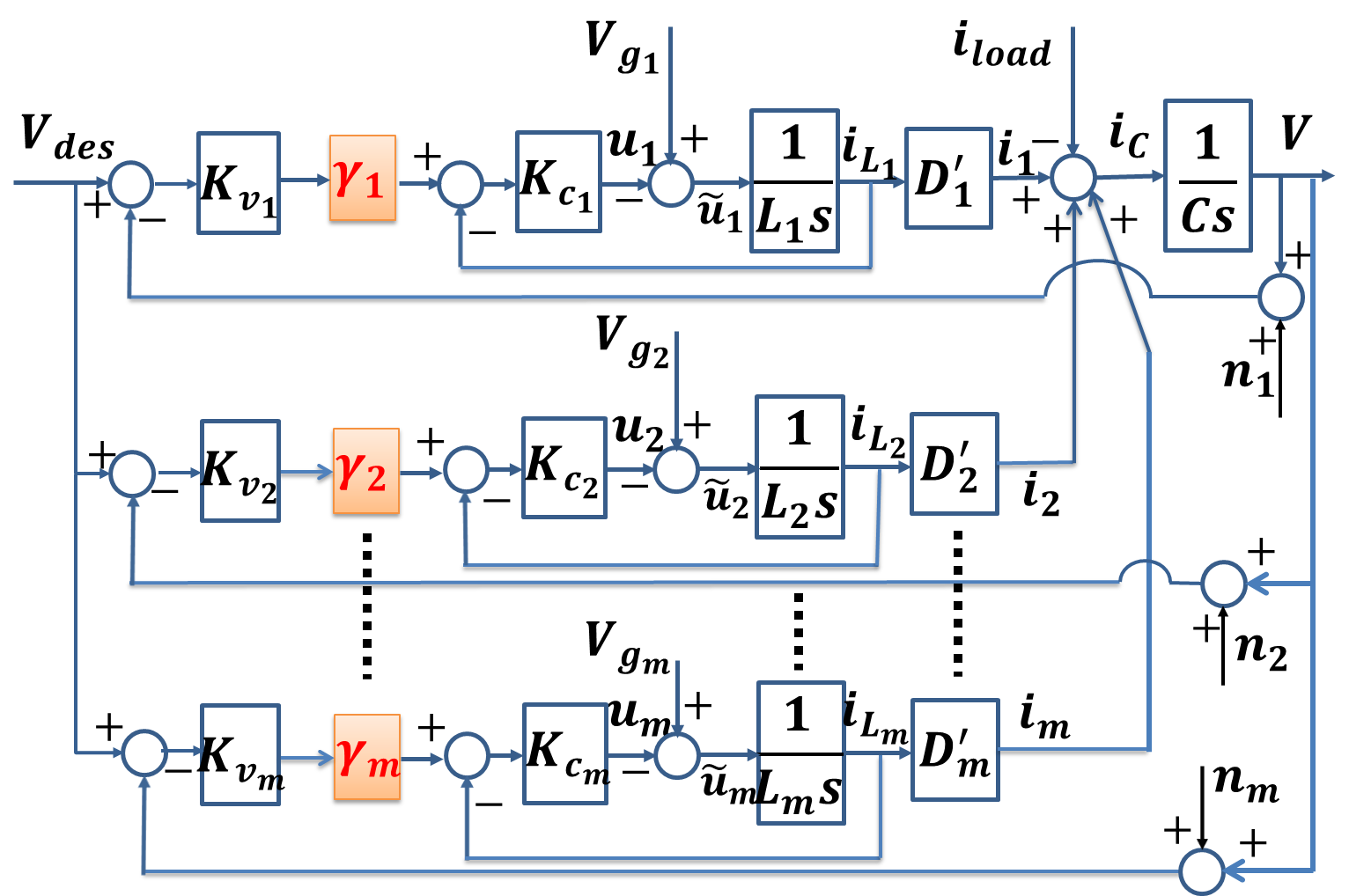}
	\vspace{-0.5em}
	\caption{{\small Control framework for a many-converters system. Note that in the proposed implementation, we adopt the same outer controller for different converters, i.e., $K_{v_1} = K_{v_2} = .... = K_{v_m} = K_v$.}}
	\label{fig:many_control}
	\vspace{-1.5em}
\end{figure}

Fig. \ref{fig:many_control} represents a decentralized inner-outer control framework for a system of $m$ parallel connected converters. Here we have incorporated a constant gain parameter $\gamma_k$ at the inner-loop of $k$th converter, the choice of which dictates power sharing as will be shown below.   After noting that the voltage-regulation objective is common to all outer controllers, in our architecture, we impose  the same outer-controller for all the converters, i.e., $K_{v_1} = K_{v_2} = .... = K_{v_m} = K_v$. This enables a significant reduction in complexity of the control design for the multi-converter system as will be shown below.

\begin{figure}[tphb]
	\centering
	\includegraphics[width=3.0in]{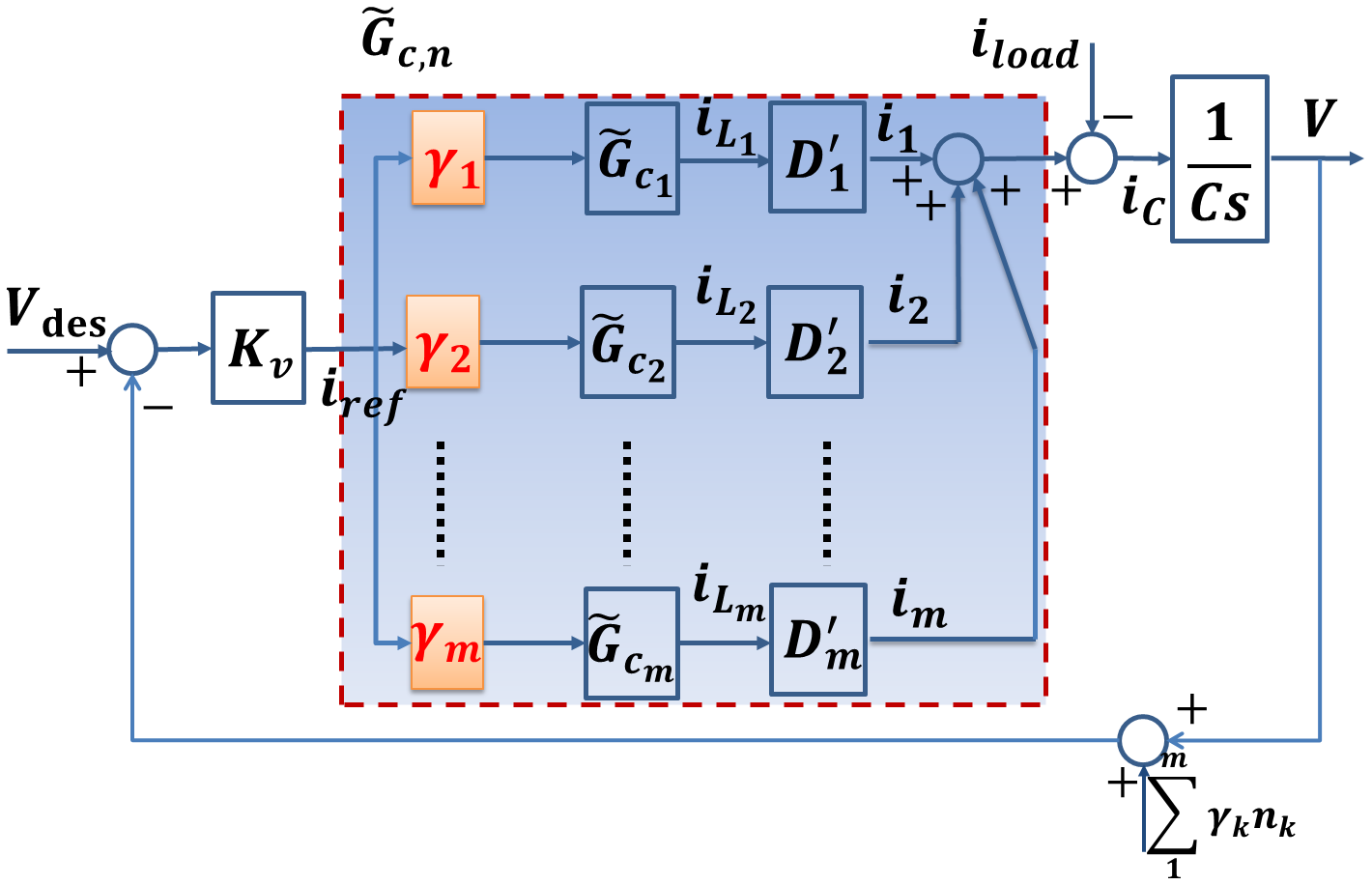}
	\vspace{-0.5em}
	\caption{{\small A multiple-converters system with shaped inner plants. Note that the shaped plants $\tilde{G}_{c_k}$ share the same denominator as $\tilde{G}_{c,n}$.}}
	\label{fig:many_control1}
	\vspace{-1.5em}
\end{figure}

First, with this assumption of $K_{v_k}=K_v$, the general decentralized architecture in Figure \ref{fig:many_control} can be simplified as in Figure \ref{fig:many_control1}. This implies that $K_v$ can be computed by solving $\mathcal{H}_\infty$-optimization problem (as discussed in the previous section) similar to the  {\em single} converter case by assuming an available  design for the summed closed inner-loop map $\tilde{G}_{c,n}$  (the {\em nominal} closed inner-loop plant)  and a {\em nominal duty-cycle} $D_n = 1-D_n'$. Then by appropriately designing individual inner-loop parameters, we can design the nominal closed inner-loop plant to be given by 
\begin{small}
\begin{equation}\label{eq:nominal_plant}
	\tilde{G}_{c,n}(s) = \left(\frac{\tilde{\omega}}{s+\tilde{\omega}}\right)\left(\frac{s^2+2\zeta_{1,n}\omega_0 s+\omega_0^2}{s^2+2\zeta_{2,n}\omega_0 s+\omega_0^2}\right),
\end{equation}
\end{small}
where the ratio $\dfrac{\zeta_{1,n}}{\zeta_{2,n}}$ determines the  tradeoff of $120$ Hz ripple between $i=\sum_ki_k$ and the capacitor current $i_C$. Note that for the cumulative closed inner-loop plant in Figure \ref{fig:many_control1} to behave as the nominal plant $\tilde{G}_{c,n}(s)$,  we require  closed inner-loop maps to sum up to the nominal closed inner-loop plant, that is $\sum_k\gamma_kD_k' \tilde G_{ck}=D_n'\tilde{G}_{c,n}$. Accordingly we design $K_{c_k}$ in each  
inner loop such  that 
\begin{small}
\begin{equation}\label{eq:nominal_plantx}
	\tilde{G}_{c_k}(s) = \left(\frac{\tilde{\omega}}{s+\tilde{\omega}}\right)\left(\frac{s^2+2\zeta_{1}^{(k)}\omega_0 s+\omega_0^2}{s^2+2\zeta_{2,n}\omega_0 s+\omega_0^2}\right),
\end{equation}
\end{small}
where $\zeta_{1}^{(k)}$ are appropriately chosen to reflect the relative tradeoff of $120$ Hz ripple among converter current outputs $i_k$. Explicit design of such $K_{c_k}$ exists and is analogous to the design in (\ref{sec:Kc}), which was obtained for the same structure of the inner-closed loop in the single-converter case.   The parameters $\gamma_k$ are designed to apportion power among the power sources, since DC gains of individual closed inner-loop plants $\gamma_k\tilde G_{ck}D'_k$ are equal to $\gamma_kD'_k$ since $\tilde G_{c_k}(j0)=1$ by design for all $k$.  We make these design specifications more precise and bring out the equivalence of the control design for the single and multiple converter systems in the following theorem.

\begin{theorem}\label{thm1}
Consider the single-converter system in Figure \ref{fig:single_inner_controller} with inductance $L$, $D'=D'_n$ and $\tilde{G_c}=\tilde{G}_{c,n}(s)$ as given 
in (\ref{eq:nominal_plant}); and  the multi-converter system described in Figures  \ref{fig:many_control} and \ref{fig:many_control1} where $\tilde{G}_{c_k}(s)$ are given by (\ref{eq:nominal_plantx}),  $\sum_k\gamma_kD_k' \tilde G_{c_k}=D_n'\tilde{G}_{c,n}$, and $\sum_k \gamma_k=1, \gamma_k>0$ for $1\leq k\leq m$. \\
\noindent 1. [Performance Equivalence]: Any outer-loop controller $K_v$ that stabilizes the single-converter system
yields identical performance when applied to the multi-converter system; more precisely, for the same exogenous inputs - the reference $V_{\text{des}}$, the load disturbance $i_{\text{load}}$, and noise $n=\sum_k\gamma_kn_k$, the steady-state regulated signals $(V_{\text{des}}-V,\dfrac{D_n'}{L}\tilde{u}, V)$ for the single-converter system are  the same as the regulated signals $(V_{\text{des}}-V,\sum_k\dfrac{D_k'}{L_k}\tilde{u}_k, V)$ for the multi-converter system.\\ 
\noindent 2. [Power Sharing]: If the parameters  $\gamma_k$, and $\zeta^{(k)}_1$, $1\leq k\leq m$ are chosen such that $\gamma_k = \frac{\alpha_k D_n'}{D_k'}\ \text{and}\ \sum\limits_{k=1}^{m}\alpha_k\zeta^{(k)}_1 = \zeta_{1,n}$, then the output current at the DC-link get divided in the ratio $\alpha_1:\alpha_2:....:\alpha_m$, where $0\leq\alpha_k\leq 1$, and $\sum\limits_{k=1}^{m}\alpha_k = 1$; more precisely the steady-state zero-frequency components $|i_1(j0)|:
|i_2(j0)|:\cdots: |i_m(j0)|$ are in the same proportion as $\alpha_1:\alpha_2:....:\alpha_m$.\\
\noindent 3. [Ripple Sharing]: If further the parameters $\zeta^{(k)}_1$ are chosen such that $\zeta^{(k)}_1 = \frac{\beta_k\zeta_{1,n}}{\alpha_k}$, where $\sum\limits_{k=1}^{m}\beta_k = 1$, and $0\leq\beta_k\leq 1, \forall k\in\{1,..,m\}$, then, the proposed design distributes the load current ripple (at 120Hz) in the ratio $\beta_1:\beta_2:...:\beta_m$; more precisely the steady-state $120$ Hz-frequency components $|i_1(j2\pi120)|:
|i_2(j2\pi120)|:\cdots:|i_m(j2\pi120)|$ are in the same proportion as $\beta_1:\beta_2:...:\beta_m$.
\end{theorem}
{ Proof: see appendix}

\section{CASE STUDIES: SIMULATIONS}\label{sec:case_studies}
In this section, we report simulation case studies, which use {\em non-ideal} components (such as diodes with non-zero breakdown voltage, IGBT switches, stray capacitances, parametric uncertainties) and switched level implementation to include nonlinearities associated with real-world experiments.
\begin{figure*}[!t]
	\begin{center}
	\begin{tabular}{ccc}
	\includegraphics[width=0.64\columnwidth]{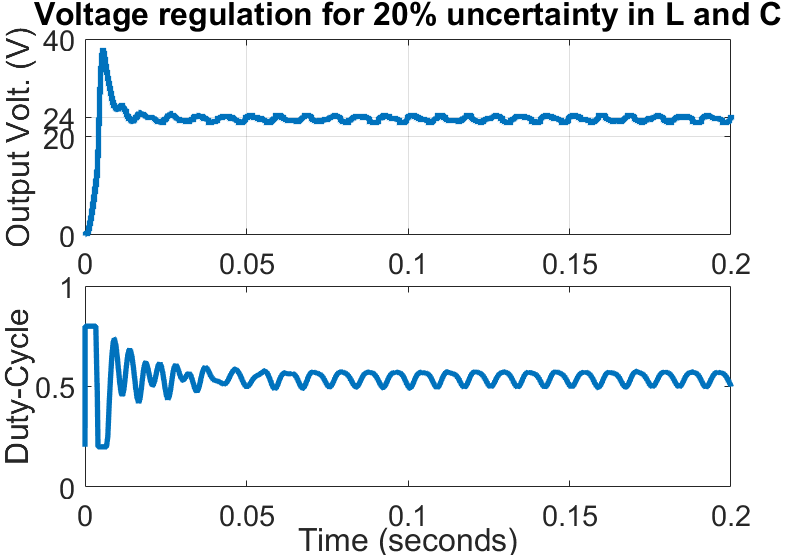}&\includegraphics[width=0.64\columnwidth]{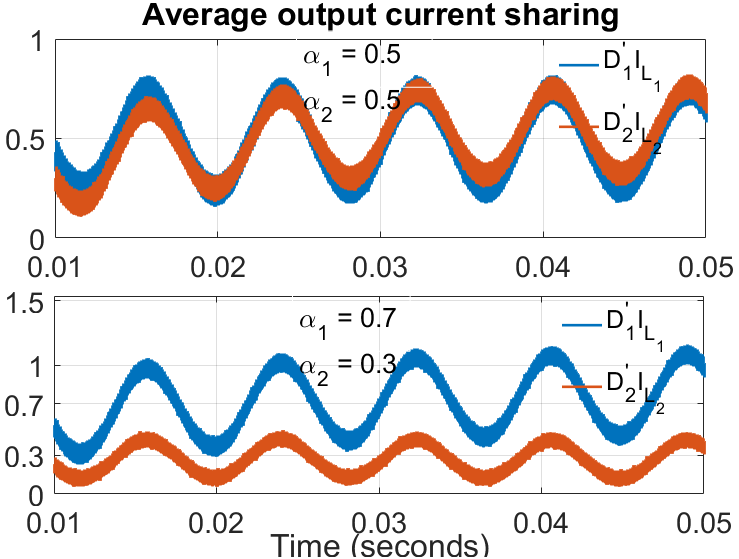}&\includegraphics[width=0.64\columnwidth]{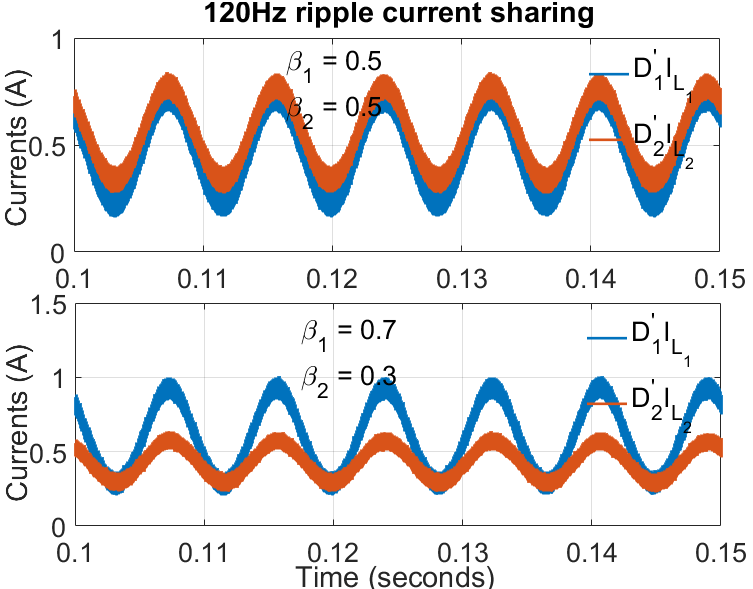}\cr
	(a)&(b)&(c)\end{tabular}
	\vspace{-0.5em}
	\caption{{\small (a)Voltage regulation in presence of modeling uncertainties and $120$ Hz ripple at the output. The inner-outer controller regulates the output voltage to the desired voltage, $V_{des} =24$V. (b) Average current (power) sharing between two converters in the ratios $1:1$ and $7:3$. (c) $\omega_0$ average ripple sharing between two converters in the ratios $1:1$ and $7:3$.}}
	\label{fig:multi_sim}
	\vspace{-2em}
	\end{center}
\end{figure*}

\subsection{Voltage regulation in presence of parametric uncertainties} \label{subsec:uncertain}
Conventional proportional-integral (PI) based control designs exhibit satisfactory performance when the actual system parameters ($L$ and $C$) lie `close' to nominal system parameters. However, a slight deviation from the nominal values of $L$ and $C$ may result in rapid degradation in the tracking performance. This issue becomes even more critical for a disturbance rejection framework, where a controller is designed without the knowledge of the uncertain load. The $\mathcal{H}_\infty$ robust control framework, where we seek an optimizing controller with guaranteed margins of robustness to modeling uncertainties will be adopted. Fig. \ref{fig:multi_sim}a shows the tracking performance of a boost converter for a $20\%$ uncertainty in both $L$ and $C$ values. The controller is designed for a boost converter with nominal $L = 2.4$mH and $C = 400\mu$F, while the actual system parameters are chosen as $L = 2$mH and $C = 500\mu$F. The input source voltage $V_g$ and the output desired voltage $V_{des}$ are chosen to be $12$V and $24$V, respectively. The design parameters for the inner-controller $K_c$ are: damping ratios $\zeta_1 = 3.2$ and $\zeta_2 = 4.5$, and $\tilde{\omega} = 2\pi 300$rad/s. The outer-controller $K_v$ is obtained by solving the stacked $\mathcal{H}_\infty$ optimization problem (see Eq. \ref{eq:stacked_h_infinity}) \cite{skogestad2007multivariable} with the weighting functions: $W_s = \frac{0.5s+2\pi 50}{s+0.06\pi 50}$, $W_u = 0.9$, and $W_t=\frac{s+2\pi 40}{0.05s+2\pi 80}$. The resulting reduced fifth-order controller $K_v$ is given by:
\begin{small}
	\begin{equation}
		K_v = \frac{0.256(s+113.9)(s+0.001)^2(s^2+4.05e4s+5.65e8)}{(s+9.56)(s^2+0.002s+4.8e-6)(s^2+9606s+8.8e7)}\nonumber
	\end{equation}
\end{small}
The load resistance is $R = 24\Omega$, and the ripple current is $I_{ripple} = 0.2\sin(2\pi 120t)$A.

\subsection{$120$ Hz ripple sharing between $i_L$ and $i_C$} \label{subsec:zeta_1_compare}
\begin{figure}[tphb]
	\centering
	\includegraphics[width=3.0in]{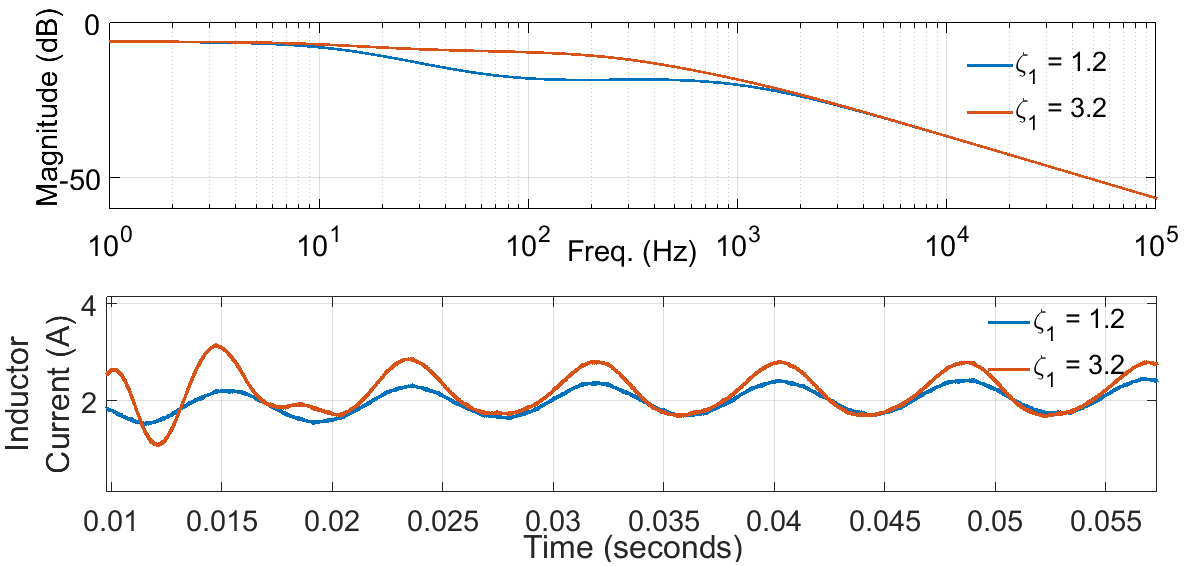}
	\vspace{-0.5em}
	\caption{{\small Bode plots of inner-shaped plants $\tilde{G_c}$ and inductor currents for different $\zeta_1$. While the average (or DC) current remains the same as desired, the smaller values of $\zeta_1$ result in relatively smaller magnitudes of $120$ Hz ripple in inductor currents.}}
	\label{fig:zeta_1_compare}
	\vspace{-2em}
\end{figure}

Fig. \ref{fig:zeta_1_compare} shows the effect of $\zeta_1$ for $120$ Hz ripple current sharing between inductor current $i_L$ and capacitor current $i_C$. Clearly, smaller values of $\zeta_1$ impart {\em notch-like} effects at $120$ Hz, thereby reducing the magnitudes of $120$ Hz ripple in inductor currents. The model and controller parameters are chosen as before.

\subsection{Average current sharing between two converters} \label{subsec:avg_sharing}
Fig. \ref{fig:multi_sim}b shows the average current sharing between two different converters with inputs $V_{g_1} = 12$V and $V_{g_2} = 10$V, respectively for two scenarios - (1) $\alpha_1 = 0.5, \alpha_2 = 0.5$, and (2) $\alpha_1 = 0.7, \alpha_2 = 0.3$. The other model and system parameters are chosen as before.

\subsection{Average $120$ Hz ripple sharing between two converters} \label{subsec:ripple_sharing}
Fig. \ref{fig:multi_sim}c shows the average $120$ Hz ripple sharing between the two converters for two scenarios - (1) $\beta_1 = 0.5, \beta_2 = 0.5$, and (2) $\beta_1 = 0.7, \beta_2 = 0.3$. The other model and system parameters are chosen as before. Note that the converters are tuned for equal average current sharing, i.e. $\alpha_1 = 0.5, \alpha_2 = 0.5$.

Thus all the objectives of the control synthesis procedure: robust voltage regulation, load power demand shared in a prescribed ratio, and the ripple current shared in a prescribed ratio are simultaneously met by our design.

\section*{APPENDIX}

\subsection{Proof of Theorem 1: Performance Equivalence}
\begin{proof}
	Let $S_n$ and $T_n$ denote the sensitivity and complementary sensitivity transfer functions of the single-converter system, respectively (as described in (\ref{eq:SnT})). For any converter $k$ in the multi-converter system in Fig. \ref{fig:many_control1}, we have $i_k = \gamma_kD_k'\tilde{G}_{c_k}K_v(V_\text{des}-V-n)$. However from Figs. (\ref{fig:many_control}) and (\ref{fig:many_control1}), we observe that $V = G_v\left(\sum\limits_{k=1}^{m}i_k-d\right)$. Thus $\sum_k\gamma_kD_k'\tilde{G}_{c_k} = D_n'\tilde{G}_{c,n}$ yields
\begin{small}
\begin{equation}\label{eq:app2}
	V_\text{des} - V = S_nV_\text{des} + G_vS_nd + T_nn, 
\end{equation}
\end{small}
which is equivalent to the map $V_\text{des} - V$ in (\ref{eq:trade_off1}) for a single-converter system. Similarly, let $G_{c_k} = \dfrac{1}{sL_k}$ denote the inner-plant in the $k$th converter, then from Fig. (\ref{fig:many_control}), it can be shown that
\begin{small}
\begin{equation}\label{eq:app3}
	%\tilde{u}_k = sL_k\gamma_k\tilde{G}_{c_k}K_v(V_\text{des}-V-n).
	\tilde{u}_k = \frac{\gamma_k}{G_{c_k}}\underset{\tilde{G}_{c_k}}{\underbrace{G_{c_k}K_{c_k}S_k}}i_{ref} = sL_k\gamma_k\tilde{G}_{c_k}K_v(V_\text{des}-V-n).
\end{equation}
\end{small}
Thus, using $\sum_k\gamma_kD_k'\tilde{G}_{c_k} = D_n'\tilde{G}_{c,n}$, we have
\begin{small}
\begin{equation}\label{eq:app4}
	\sum_k\frac{D_k'}{L_k}\tilde{u}_k = sD_n'\tilde{G}_{c,n}K_v(V_\text{des}-V-n) = \frac{D_n'}{L}\tilde{u},
\end{equation}
\end{small}
which establishes the required equivalence.
\end{proof}

\subsection{Proof of Theorem 1: Power Sharing}
\begin{proof}
	Note that from Fig. (\ref{fig:many_control1}), we have
	\begin{small}\begin{equation}\label{eq:app5}
		i_k(s) = \gamma_kD_k'\tilde{G}_{c_k}(s)i_\text{ref}(s).
	\end{equation}\end{small}
	Thus, using the fact that $|\tilde{G}_{c_k}(j0)| = 1$ and with the given choice of the parameter $\gamma_k = \left(\alpha_kD_n'/D_k'\right)$, we obtain $\left(|i_k(j0)|/\sum_k|i_k(j0)|\right) = \left(\gamma_kD_k'/\sum_k\gamma_kD_k'\right) = \alpha_k$. Thus, the steady-state zero-frequency component of the output current at the DC-link gets divided in the ratio $\alpha_1:\alpha_2:....:\alpha_m$.
\end{proof}

\subsection{Proof of Theorem 1: $120$ Hz ripple sharing}
\begin{proof}
	From (\ref{eq:app5}) and observing that $|\tilde{G}_{c_k}(j\omega_0)| = \left|\dfrac{\tilde{\omega}}{j\omega_0 + \tilde{\omega}}\right|\dfrac{\zeta_1^{(k)}}{\zeta_{2,n}}$, the ratio of $120$ Hz ripple magnitude in steady-state is given by $\left(|i_k(j\omega_0)|/\sum\limits_{k'=1}^{m}|i_{k'}(j\omega_0)|\right) = \left(\gamma_kD_k'\zeta^{(k)}_1/\sum\limits_{k'=1}^{m}\gamma_{k'}D_{k'}'\zeta^{(k')}_1\right)$. Substituting $\gamma_kD_k$ $= \alpha_kD_n'$ and $\sum\limits_{k=1}^{m}\alpha_k\zeta^{(k)}_1$ $= \zeta_{1,n}$ yields $\left(|i_k(j\omega_0)|/\sum\limits_{k'=1}^{m}|i_{k'}(j\omega_0)|\right)$ $= \left(\alpha_k\zeta^{(k)}_1/\zeta_{1,n}\right)$. But, by our choice of the damping parameters, $\zeta^{(k)}_1 = \left(\beta_k\zeta_{1,n}/\alpha_k\right)$, yields $\left(|i_k(j\omega_0)|/\sum\limits_{k'=1}^{m}|i_{k'}(j\omega_0)|\right) = \beta_k$. Thus, the ripple currents get divided in the ratios, $\beta_1:\beta_2:...:\beta_m$.
\end{proof}

\bibliographystyle{IEEEtran} 
\bibliography{myRef}

%%%%%%%%%%%%%%%%%%%%%%%%%%%%%%%%%%%%%%%%%%%%%%%%%%%%%%%%%%%%%%%%%%%%%%%%%%%%%%%%

\end{document}